\documentclass{ifacconf}
\usepackage{graphicx}    
\usepackage{natbib}

\usepackage{pgfplots}
\pgfplotsset{compat=1.18} 

\usepackage{wrapfig}   
\usepackage{amsmath,amssymb,mathrsfs,bbold,stmaryrd}

\usepackage{booktabs}
\usepackage{algorithm}
\usepackage{comment}
\newtheorem{definition}{Definition}
\newtheorem{proposition}{Proposition}
\newtheorem{theorem}{Theorem}

\newtheorem{example}{Example}

\newcommand{\Sr}{{{\rm Sres}}}
\newcommand{\sr}{{{\rm sres}}}

\newcommand{\Sth}{{{\rm StHa}}}
\newcommand{\sth}{{{\rm stha}}}
\newcommand{\Res}{{{\rm Res}}}
\newcommand{\Lc}{{{\rm Lc}}}

\usepackage{xcolor}

\begin{document}
\begin{frontmatter}

\title{Symbolic and Numerical Tools for $L_{\infty}$-Norm Calculation} 

\author[SUAD]{Grace Younes}
\author[INRIA]{Alban Quadrat} 
\author[INRIA]{Fabrice Rouillier}
\address[SUAD]{Sorbonne University Abu Dhabi, SAFIR, Abu Dhabi, UAE \\ 
\texttt{grace.younes@sorbonne.ae}}

\address[INRIA]{Sorbonne Universit\'e and Universit\'e de Paris, CNRS, IMJ-PRG, Inria Paris, F-75005 Paris, France \\ 
\texttt{\{alban.quadrat, fabrice.rouillier\}@inria.fr}}

\begin{abstract}
The computation of the \( L_\infty \)-norm is an important issue in $H_{\infty}$ control, particularly for analyzing system stability and robustness. This paper focuses on symbolic computation methods for determining the \( L_{\infty} \)-norm of finite-dimensional linear systems, highlighting their advantages in achieving exact solutions where numerical methods often encounter limitations. Key techniques such as Sturm-Habicht sequences, Rational Univariate Representations (RUR), and Cylindrical Algebraic Decomposition (CAD) are surveyed, with an emphasis on their theoretical foundations, practical implementations, and specific applicability to \( L_{\infty} \)-norm computation.
A comparative analysis is conducted between symbolic and conventional numerical approaches, underscoring scenarios in which symbolic computation provides superior accuracy, particularly in parametric cases. Benchmark evaluations reveal the strengths and limitations of both approaches, offering insights into the trade-offs involved. Finally, the discussion addresses the challenges of symbolic computation and explores future opportunities for its integration into control theory, particularly for robust and stable system analysis.
\end{abstract}

\end{frontmatter}
\section{Introduction}
\label{sec:introduction}

The computation of the \( L_{\infty} \)-norm is an important problem in $H_{\infty}$ control, with key applications in stability analysis, robustness evaluation, and system optimization. This norm provides valuable insights into the worst-case behavior of a system, making its accurate computation essential. However, despite its importance, the computation of the \( L_{\infty} \)-norm remains challenging, particularly for parameterized systems or when high precision is required.

Existing numerical methods, such as bisection algorithms and optimization-based approaches, often struggle with these challenges. Specifically, they may fail to provide reliable results due to issues such as limited precision, ill-conditioning, or the inability to identify all real roots of the governing polynomial equations. These limitations underscore the need for alternative approaches that can deliver exact and reliable solutions, even in complex settings.

Symbolic computation offers a promising solution to these challenges. By leveraging algebraic techniques, symbolic methods can compute the \( L_{\infty} \)-norm with high precision and provide exact solutions, particularly in parametric cases where numerical methods often fail. Nevertheless, symbolic computation is not without its drawbacks, as it is often computationally intensive and struggles to scale for larger systems. Overcoming these limitations is crucial for fully realizing the potential of symbolic methods in practical control applications. 

This paper provides a detailed exploration of symbolic computation tools and methods for \( L_{\infty} \)-norm computation, emphasizing their theoretical foundations, practical implementations, and relevance to control theory. Specifically, we present an overview of key symbolic tools, including \emph{Sturm-Habicht sequences}, \emph{Rational Univariate Representations} (RUR), and \emph{Cylindrical Algebraic Decomposition} (CAD), with a focus on their applicability to \( L_{\infty} \)-norm computation. Furthermore, we compare symbolic methods with existing numerical approaches, highlighting scenarios in which symbolic computation excels, particularly in handling parametric systems. We also discuss the limitations of symbolic computation and identify potential research directions for integrating symbolic techniques into practical control system analysis.

The paper is organized as follows. Section~\ref{sec:infinity-norm} introduces the \( L_{\infty} \)-norm problem in control theory, covering its definition, numerical methods, and challenges. Section~\ref{sec:symbolic_tools} discusses symbolic and hybrid symbolic-numeric approaches for \( L_{\infty} \)-norm calculation. Section~\ref{sec:versus} compares symbolic and numerical methods through examples and benchmarks that evaluate the performance of numerical methods in {\tt Matlab} and {\tt Maple} versus symbolic methods in {\tt Maple}. Finally, Section~\ref{sec:conclusion} summarizes the findings and suggests future research directions.

\section{The $L_{\infty}$-Norm in Control Theory}
\label{sec:infinity-norm}

We first introduce standard notations and definitions. Let \(\mathbb{k}\) be an arbitrary field, typically chosen as \(\mathbb{Q}\) or \(\mathbb{R}\). For a polynomial \(P \in \mathbb{k}[x, y]\), we denote by \(\operatorname{Lc}_{\mathrm{var}}(P)\) the \emph{leading coefficient} of \(P\) with respect to the variable \(\mathrm{var} \in \{x, y\}\), and by \(\deg_{\mathrm{var}}(P)\), the \emph{degree} of \(P\) in the variable \(\mathrm{var} \in \{x, y\}\). The \emph{total degree} of \(P\) is denoted by \(\deg(P)\). 
Next, consider the projection maps from the real plane \(\mathbb{R}^2\) onto its coordinate axes. The map \(\pi_x: \mathbb{R}^2 \to \mathbb{R}\) is the projection onto the \(x\)-axis, defined by \(\pi_x(a, b) = a\) for all \((a, b) \in \mathbb{R}^2\). Similarly, the \(y\)-projection \(\pi_y: \mathbb{R}^2 \to \mathbb{R}\) maps \((a, b)\) to \(b\). In this context, we say that \(a\) is the \(x\)-projection of the point \((a, b)\), and \(b\) is its \(y\)-projection.
For polynomials \(P, Q \in \mathbb{k}[x, y]\), let \({\rm gcd}(P, Q)\) denote their greatest common divisor. The ideal of \(\mathbb{k}[x, y]\) generated by \(P\) and \(Q\) is denoted by \(I := \langle P, Q \rangle\). The affine algebraic set defined by \(I\) over a field \(\mathbb{K}\) (where \(\mathbb{K}\) is a field containing \(\mathbb{k}\)) is given by $
V_{\mathbb{K}}(I) := \{(x, y) \in \mathbb{K}^2 \; | \; \forall \, R \in I, \; R(x, y) = 0\}. $
Finally, let \(\mathbb{C}_+ := \{s \in \mathbb{C} \; | \; \Re(s) > 0\}\) represent the \emph{open right half-plane} of the complex plane \(\mathbb{C}\).  Additionally, for any matrix \(A\), \(\det(A)\) refers to the \emph{determinant} of \(A\).

\subsection{Characterization of the $L_{\infty}$-Norm}
\begin{definition}[\cite{Doyle,Zhou}]
Let ${RH}_{\infty}$ be the $\mathbb{R}$-algebra of all the proper and stable rational functions with real coefficients,
namely:
\begin{align*}
   {RH}_{\infty}:= \left \{\dfrac{n}{d} \; | \; n, d \in \mathbb{R}[s], \; {\rm gcd}(n, d)=1, \right.\\ \left.\deg_s(n) \leq \deg_s(d), V_{\mathbb{C}}(\langle d \rangle) \cap \mathbb{C}_+=\emptyset \right \}.
\end{align*}  

\end{definition}

Each $g \in {RH}_{\infty}$ is holomorphic and bounded on $\mathbb{C}_+$ (i.e. \(\|g\|_\infty := \sup_{s\in\mathbb{C}_+}|g(s)|<\infty\)), 
and the maximum modulus principle yields
$
\|g\|_\infty = \sup_{\omega\in\mathbb{R}} |g(i\omega)|.
$
Thus, the restriction \( g|_{i\mathbb{R}} \) belongs to the \emph{Lebesgue space} \(L_\infty(i\mathbb{R})\), and more precisely to the algebra ${RL}_{\infty}$ of 
real rational functions on the imaginary axis $i \, \mathbb{R}$ which are proper and have no poles on $i \, \mathbb{R}$, or simply, the algebra of real rational functions with no poles on $i \, \mathbb{P}^1(\mathbb{R})$, where 
$\mathbb{P}^1(\mathbb{R}):=\mathbb{R} \cup {\infty}$.

We extend the $L_\infty$-norm from functions in $RH_\infty$ (resp., ${RL}_{\infty}$) to matrices. For $G \in {RH}_{\infty}^{u \times v}$ (resp., $G \in {RL}_{\infty}^{u \times v}$, $\mathbb{R}(s)^{u \times v}$), 
let $\Bar{\sigma}\left ( \cdot \right)$ denote the largest \emph{singular value of a complex matrix}. Then, we define:
$$\parallel G \parallel_{\infty} \, := \,  \sup_{s \in \mathbb{C}_+} \Bar{\sigma}\left(G(s) \right) \;  \big(
\mbox{resp.,} \, \parallel G \parallel_{\infty} \, := \,  \sup_{\omega \in \mathbb{R}} \Bar{\sigma}\left(G(i \, \omega) \right) 
\big).$$ 
For $G \in {RH}_{\infty}^{u \times v}$, this yields $\parallel G \parallel_{\infty} \, = \, \sup_{\omega \in \mathbb{R}} \Bar{\sigma}\left(G(i \, \omega) \right)$. 
The \emph{conjugate} of \(G \in \mathbb{R}(s)^{u\times v}\) is defined as
\[
\tilde{G}(s) := G^T(-s).
\]

\begin{proposition}[\cite{kanno2006validated}]\label{prop:Kanno-Smith}
Let $G\in \mathbb{R}(s)^{u\times v}$ with $G|_{i\mathbb{R}}\in RL_\infty^{u\times v}$ and $\gamma>0$. Define   
$$\Phi_{\gamma}(s)= \gamma^2 \, I_v-\Tilde{G}(s) \, G(s).$$ Then, $\gamma > \, \parallel G \parallel_{\infty}$ if and only if $\gamma
> \Bar{\sigma}\left(G(i \, \infty)\right)$ and $\det(\Phi_{\gamma}(i \, \omega)) \neq 0$ for all $\omega \in \mathbb{R}$.    
\end{proposition}
\subsection{Overview of Numerical Methods}

Unlike the $L_2$-norm, no closed-form expression exists for the $L_{\infty}$-norm of finite-dimensional systems (\cite{Doyle,Zhou}). Thus, its computation relies on efficient numerical methods such as \emph{bisection} and \emph{Hamiltonian eigenvalue computations} (\cite{Boyd,Bruinsma}). In \cite{Boyd}, the authors related the singular values of the transfer matrix on the imaginary axis to the eigenvalues of a Hamiltonian matrix, leading to a bisection algorithm for the \(L_{\infty}\)-norm. Similarly, in \cite{Bruinsma}, the authors proposed a fast algorithm leveraging this relationship.
Recent methods focus on localizing real roots of bivariate polynomials. For instance, \cite{kanno2006validated} used Sturm chain tests for validation, while \cite{Belur} applied structured linearization of \emph{Bezoutian matrices} to compute isolated common zeros.

\subsection{Challenges and Limitations of Numerical Methods}

Although efficient, numerical methods assume fixed coefficients, and thus, struggle with parameterized transfer matrices (\cite{chen2013computing}).
Additionally, floating-point arithmetic can cause precision issues in ill-conditioned systems or near singularities, where small errors may destabilize results (\cite{kanno2006validated}).
Iterative approaches (e.g., bisection, eigenvalue methods) may suffer from convergence issues, often relying on initial guesses and converging to incorrect local optima (\cite{Boyd,Bruinsma}).
These challenges call for methods that guarantee accuracy in parameterized or ill-conditioned problems. Symbolic computation offers a promising solution with its exact, certified approaches.

\section{Symbolic Computation Tools}
\label{sec:symbolic_tools}

\subsection{Symbolic Reformulation of the $L_{\infty}$-Norm Computation}

Since \(\det(\Phi_\gamma(i \, \omega))\) is a real rational function (in fact, of \(\omega^2\) and \(\gamma^2\)), we write it as
\(
\det(\Phi_\gamma(i\omega))=\frac{n(\omega,\gamma)}{d(\omega)},
\)
with \(n\in\mathbb{R}[\omega,\gamma]\) and \(d\in\mathbb{R}[\omega]\) coprime. Thus, the maximal singular value of \(G(i \, \omega)\) is the largest real \(\gamma\) for which there exists some real \(\omega\) such that \(n(\omega,\gamma)=0\).

Since \(G\) has no poles on the imaginary axis, \(d(\omega)\neq 0\) for all \(\omega\in\mathbb{R}\). Thus, the \(L_{\infty}\)-norm of \(G\) is given by the maximal \(\gamma\) for which \(n(\omega,\gamma)=0\) has a real solution. \emph{Hence, we study the \(\gamma\)-extremal (critical) points of the curve} \(\mathcal{C}=\{(\omega,\gamma)\in\mathbb{R}^2 \mid n(\omega,\gamma)=0\}\). A \(\gamma\)-critical point occurs when \(n(\omega,\gamma)\) reaches a local maximum, minimum, or a turning point, i.e., when \(n(\omega,\gamma)\)  and its derivative \(\partial n/\partial\omega\) share a common factor.

{\begin{definition}\label{def:curve_n}
    Let $\gamma > 0$ and $G \in \mathbb{R}(s)^{u \times v}$ be such that $G_{|i \, \mathbb{R}} \in {RL}_{\infty}^{u \times v}$. The numerator $n \in \mathbb R [\omega,\gamma]$ of $\det\left(\gamma^2 \, I_v-G^T(-i \, \omega) \, G(i \, \omega)\right)$
        is called the \emph{polynomial associated with the $L_{\infty}$-norm}.
\end{definition}
}

\begin{proposition}[\cite{bouzidi2021computation}]\label{prop:caracterization_infty_norm}
Let $G \in {RL}_{\infty}^{u \times v}$ and $n \in \mathbb{R}[\omega, \gamma]$ be {the polynomial associated with the $L_{\infty}$-norm}. We denote by $\Bar{n} \in \mathbb{R}[\omega,\gamma]$ the square-free part of $n$, that is, the polynomial obtained by removing any repeated factors from \( n \).
Then, we have: 
$$\parallel G \parallel_{\infty} =\max \left \{ \pi_{\gamma} \left ( V_{\mathbb{R}} \left(
\left \langle 
\Bar{n}, \frac{\partial \Bar{n}}{\partial \omega} 
\right \rangle 
\right)\right) \cup V_{\mathbb{R}}\left(\langle \operatorname{Lc}_{\omega}(\Bar{n}) \rangle \right)  \right \},$$
where the real roots of \( \operatorname{Lc}_{\omega}(\Bar{n}) \) correspond to values of \( \gamma \) where the real curve \( \Bar{n} = 0 \) has horizontal asymptotes.

\end{proposition}

Assume \(n(\omega,\gamma)\) is square-free. Then, by Proposition~\ref{prop:caracterization_infty_norm}, for \(G\in RL_\infty^{u\times v}\), computing \(\|G\|_\infty\) reduces to finding the maximal \(\gamma\)-coordinate among the real solutions of the system \(\Sigma=\{n(\omega,\gamma)=0,\; \partial n(\omega,\gamma)/\partial \omega=0\}\).
This amounts to finding the largest positive \(\gamma_\star\) such that \(n(\omega,\gamma_\star)=0\) has a real solution.

Using symbolic-numeric techniques,  this problem was studied in \cite{bouzidi2021computation, kanno2006validated} for \(n\in\mathbb{Q}[\omega,\gamma]\) to compute the \(L_{\infty}\)-norm. In \cite{quadrat:hal-04646145}, the approach was extended to \(n\in\mathbb{Z}[\alpha][\omega,\gamma]\), with \(\alpha=(\alpha_1,\ldots,\alpha_d)\) denoting real parameters. The symbolic tools employed are detailed in the following sections.

For a transfer matrix \( G \in \mathbb{R}(s)^{p \times m} \), the \( L_{\infty} \)-norm satisfies \( \| G \|_{\infty} < \gamma \) if and only if \( \Phi_{\gamma}(i \, \infty) > 0 \) and \( \Phi_{\gamma}(i \, \omega) \) is non-singular for all \( \omega \in \mathbb{R} \), where \( \Phi_{\gamma}(s) = \gamma^2 I_m - G^T(-s) \, G(s) \). 
If $G = C (s I_n - A)^{-1} B + D$, where $A$ has no purely imaginary eigenvalues (since, otherwise, \(\| G \|_{\infty}=\infty\)), considering a minimal  
realization of $\Phi_{\gamma}^{-1}$  as 
$$\Phi_{\gamma}^{-1}=E \, (s \, I-H_{\gamma})^{-1} \, F+G,$$
where \( H_{\gamma} \) is a Hamiltonian matrix depending of $A, \, B, \, C$, $D$ and $\gamma$, then the poles of $\Phi_{\gamma}^{-1}(i \, \omega)$, i.e., the solutions of $n(\omega, \gamma)=0$, are then the eigenvalues of $H_{\gamma}$ on the imaginary axis. This remark yields the following result  (\cite{Boyd,Bruinsma,Zhou}).

\begin{theorem}\label{thm:Hamiltonian_imaginary_eigenvalues}
Let \( G = C (s I_n - A)^{-1} B + D \), where \( A \) has no eigenvalues on the imaginary axis, i.e., \( G \in RL_{\infty}^{p \times m} \), and let \( \gamma > 0 \) with \( R = \gamma^2 I_m - D^T D \). Then, \( \| G \|_{\infty} < \gamma \) is equivalent to  \(\overline{\sigma}(D) < \gamma \) and \( H_{\gamma} \) has no imaginary eigenvalues.
\end{theorem}

This result reduces \(L_{\infty}\)-norm computation to analyzing the eigenvalues of \(H_\gamma\), bridging numerical methods and our symbolic approach.

\subsection{Key Symbolic Computation Tools}\label{symbolic-tools}

Symbolic methods are crucial for the reformulated $L_{\infty}$-norm problem, especially for parameterized systems or when certified results are needed. This subsection introduces key symbolic methods in polynomial and system analysis, emphasizing theory and computation.

\textbf{Greatest Common Divisor (GCD).}
The \emph{greatest common divisor} (GCD) of two polynomials is a fundamental concept in symbolic computation, enabling the simplification of polynomial systems, root isolation, and multivariate-to-univariate reductions. Formally, let \(\mathbb{D}\) be an \emph{integral domain} $-$ a commutative ring with no non-trivial zero-divisors $-$ typically \(\mathbb{Z}[x]\), \(\mathbb{Z}[x, y]\), or \(\mathbb{Z}\). The GCD of two polynomials \(P\) and \(Q\) in \(\mathbb{D}[x]\), denoted \(\gcd(P, Q)\), is the polynomial of highest degree that divides both \(P\) and \(Q\) and is unique up to a unit in \(\mathbb{D}\) (\cite{cox2013ideals}).
For Euclidean domains (integral domains permitting division with remainder) such as \(\mathbb{Z}[x]\) or \(\mathbb{k}[x]\), the Euclidean Algorithm computes the GCD via successive divisions (\cite{basu2006existential}). This cornerstone of computational algebra is both efficient and structure-preserving.

\textbf{Resultants.}

The \emph{resultant} is key concept in symbolic computation for eliminating variables in multivariate systems.

For \(P, \, Q \in \mathbb{D}[x]\) of degrees  \(p\) and \(q\), the resultant \(\operatorname{Res}(P,Q,x)\) is the determinant of the \((p+q)\times(p+q)\) \emph{Sylvester matrix} constructed from their coefficients. It vanishes if and only if \(P\) and \(Q\) share a nontrivial factor (\cite{cox1998d}). Resultants relate to \emph{discriminants}. For a univariate polynomial \(P\), its \emph{discriminant} is
\[
\operatorname{discrim}(P) = (-1)^{\frac{p(p-1)}{2}}\operatorname{Lc}(P)^{-1}\operatorname{Res}(P,P',x),
\]
where \(P'\) is the derivative with respect to $x$. It vanishes if and only if \(P\) has multiple roots (\cite{Basu}).

\begin{example}\label{ex:gcd}
Consider \( P = x^4 - 1 \) and \( Q = x^6 - 1 \). The Sylvester matrix's determinant evaluates to \( 0 \), confirming a common factor. Using the Euclidean algorithm, the greatest common divisor of $P$ and $Q$ is \( x^2 - 1 \).    
\end{example}

Let \(P, \, Q \in \mathbb{Z}[x,y]\) be two bivariate polynomials, viewed as univariate in \(x\) with degrees \(p\) and \(q\) with coefficients in \(\mathbb{Z}[y]\). The resultant \(\Res(P,Q,x)=\det(L)\in\mathbb{Z}[y]\), where \(L\) is the \((p+q) \times (p+q)\) Sylvester matrix of \(P\) and \(Q\), captures the \(y\)-coordinates of the solutions of \(\{P=0, \, Q=0\}\).

\begin{proposition}\label{projection poperty of resultant polynomial}[\cite{cox1998d}]
Let \( P, \, Q \in \mathbb{Z}[x, y] \) and \(\Res(P, Q, v)\) be their resultant with respect to \(v \in \{x, y\}\). Suppose \(\Lc_v(P) = a_p\) and \(\Lc_v(Q) = b_q\), with \(a_p, b_q \in \mathbb{Z}[u]\), \(u \in \{x, y\} \setminus \{v\}\). For \(\alpha \in \mathbb{C}\), the following assertions are equivalent:
\begin{enumerate}
    \item \(\Res(P, Q, v)(\alpha) = 0\).
    \item \(a_p(\alpha) = b_q(\alpha) = 0\), or there exists \(\beta \in \mathbb{C}\) such that \((\alpha, \beta) \in V_{\mathbb{C}}\big(\langle P, Q \rangle\big)\).
\end{enumerate}
\end{proposition}
\begin{example}
Let \(P=x \, y-1\) and \(Q=x^2+y^2-4\). Viewed as univariate in \(x\) (with coefficients in \(\mathbb{Q}[y]\)), then
\[
\operatorname{Res}(P,Q,x)=\det\begin{pmatrix} y & -1 & 0 \\ 0 & y & -1 \\ 1 & 0 & y^2-4 \end{pmatrix}=y^4-4 \, y^2+1.
\]
The roots of \(y^4-4 \, y^2+1=0\) give the \(y\)-projections of the solutions of the polynomial system \(\{P=0,\,Q=0\}\).
\end{example}
\begin{proposition}\label{Prop:Res-norm}
Let \(G\in RL_\infty^{u\times v}\) and \(n\in\mathbb{R}[\omega,\gamma]\) be the polynomial associated with the \(L_{\infty}\)-norm. Then, as per Propositions~\ref{prop:caracterization_infty_norm} and \ref{projection poperty of resultant polynomial}, if \(n\) is viewed as a univariate polynomial in \(\omega\), then 
\(\|G\|_\infty\)
is the maximal real root \(\gamma_\star\) of 
\(
\operatorname{Res}\Bigl(n,\frac{\partial n}{\partial \omega},\omega\Bigr)
\)
for which \(n(\omega,\gamma_\star)\) has a real root or \(\operatorname{Lc}_\omega(n)(\gamma_\star)=0\).
\end{proposition}

\begin{example}\label{resultants_and_norm}
Consider \(G=\frac{1}{2 \, s^2+3 \, s+2}\). As per Definition~\ref{def:curve_n},
\(
n=4 \, \gamma^2 \, \omega^4+\gamma^2 \, \omega^2+4 \, \gamma^2-1
\)
and 
\(
R=\operatorname{Res}(n,\partial n/\partial \omega,\omega)
\)
has square-free part
\(
\gamma \, (2 \, \gamma-1)(2 \, \gamma+1)(63 \, \gamma^2-16),
\)
with real roots 
\(
\left\{-\sqrt{\frac{16}{63}}, -\frac{1}{2}, 0, \frac{1}{2}, \sqrt{\frac{16}{63}}\right\}.
\)
Here, \(0\) (the root of \(\operatorname{Lc}_\omega(n)\)) is an asymptote. Since 
\(
n\Bigl(\omega,\sqrt{\frac{16}{63}}\Bigr)=\frac{1}{63} \, (8 \, \omega^2+1)^2
\)
has no real roots and 
\(
n\Bigl(\omega,\frac{1}{2}\Bigr)=\frac{1}{4} \, \omega^2(4 \, \omega^2+1)
\)
has a real root (\(\omega=0\)), we conclude \(\|G\|_{\infty}=\frac{1}{2}\).
\end{example}

Example~\ref{resultants_and_norm} raises the question: How do we count the real roots of \( n(\omega,\gamma_\star) \) when \(\gamma_\star\) is an algebraic number? In this context, \emph{subresultants} and \emph{Sturm-Habicht sequences} are essential as the variation in their leading coefficients' signs directly indicates the number of real roots.

\textbf{Subresultants and Sturm-Habicht Sequences.} The \emph{subresultant sequence} of two univariate polynomials \(P, \, Q \in \mathbb{D}[x]\), denoted by $\{\Sr_i(P, Q, x)\}$, generalizes the Euclidean algorithm for these polynomials. Suppose that $\operatorname{deg}(P) = p$, $\operatorname{deg}(Q) = q$ and, without loss of generality, {$p<q$. Let $\lambda = \min(p,\ q)$}. Then, $\operatorname{Sres}_{\lambda}(P, Q, x)=P,\; \operatorname{Sres}_{\lambda+1}(P, Q, x)=Q$, and, for $0 \leq i < \lambda$, each subresultant $\operatorname{Sres}_{i}$ is a polynomial in $\mathbb{D}[x]$ derived as a remainder of a modified Euclidean algorithm and of coefficients determined by the determinants of structured submatrices of the Sylvester matrix of \(P\) and \(Q\) (\cite{basu2006existential,kahoui2003elementary}). We denote by $\operatorname{sres}_i(P,Q,x)$ the leading coefficient of $\operatorname{Sres}_i(P,Q,x)$.

One can prove that $\operatorname{Sres}_0(P, Q, x)= \operatorname{sres}_0(P, Q, x) = \operatorname{Res}(P, Q, x) \in \mathbb{D}$. Furthermore, if \(\sr_i(P, Q, x) = 0\) for \(i = 0, \ldots, k-1\) and \(\sr_k(P, Q, x) \neq 0\), then we have \(\gcd(P, Q) = \Sr_k(P, Q, x)\) over the quotient field \(\mathbb{F}_{\mathbb{D}}\) of $\mathbb{D}$ (\cite{basu2006existential,kahoui2003elementary}). 

\begin{example}
Consider again \(P = x^4 - 1\) and \(Q = x^6 - 1\) as in Example~\ref{ex:gcd}.  By computing the subresultant sequence of \(P\) and \(Q\) with respect to the variable \(x\), we get \(\Sr_0(P, Q, x) = 0\), \(\Sr_1(P, Q, x) = 0\), \(\Sr_2(P, Q, x) = x^2 - 1\), and \(\Sr_3(P, Q, x) = -x^2 + 1\).
Hence, the least non-vanishing subresultant is \(\Sr_2(P, Q, x) = x^2 - 1\), confirming that \(\gcd(P, Q) = x^2 - 1\).
\end{example}

Furthermore, the specialization property of subresultants ensures that for a \emph{ring homomorphism} \(\phi: \mathbb{D} \to \mathbb{D}'\), we have \(\phi(\Sr_j(P, Q, x)) = \Sr_j(\phi(P), \phi(Q), x)\). This property, established in \cite{gonzalez1994specialisation, gonzalez1990specialisation, kahoui2003elementary}, guarantees that subresultants remain stable under specialization as illustrated in the next example.

\begin{example}\label{specialization example}

Let $$P(x, y) = x^4 - (y + 2)\, x^3 + (2\, y + 1)\, x^2 - (y + 2)\, x + 2\, y,$$ $Q(x, y) =\partial P (x, y)/\partial x$,
and $\phi: \mathbb Z[x, y] \longrightarrow \mathbb Z[x, y]$ be defined by {$\phi(P(x, y))= P(x, {\alpha})$}, where ${\alpha} \in \mathbb Q$.
Then, 
\begin{itemize}
    \item $\Sr_0(\phi(P), \phi(Q), x) = \Sr_0(P, Q, x)(x, \alpha) \vspace{1mm}\\ 
    \hspace*{3.1cm} =  -100\, (\alpha - 2)^2\, (\alpha^2 + 1)^2 $, 
    \item  $\Sr_1(\phi(P), \phi(Q), x) = \Sr_1(P, Q, x)(x, \alpha) \vspace{1mm}\\ = (-2\, \alpha^4 + 40\, \alpha^3 - 74\, \alpha^2 + 40\, \alpha + 128)\, x  \vspace{1mm}\\
    \hspace*{4mm} -4\, (\alpha + 2)\, (4\, \alpha^3 - 13\, \alpha^2 + 24\, \alpha - 3) $, 
    \item $\Sr_2(\phi(P), \phi(Q), x) = \Sr_2(P, Q, x)(x, \alpha) \vspace{1mm}\\ = (-3 \, \alpha^2 + 4 \, \alpha - 4)\, x^2 + 2\, (\alpha + 2)\, (2\, \alpha - 5)\, x \vspace{1mm} \\
    \hspace*{4mm} - \alpha^2 + 28\, \alpha - 4  $, 
    \item $\Sr_3(\phi(P), \phi(Q), x) = \Sr_3(P, Q, x)(x, \alpha) = Q(x, \alpha) \vspace{1mm} \\ =  4\, x^3 - (3\, \alpha + 6)\, x^2 + (4\, \alpha + 2)\, x - \alpha - 2 $,
    \item $\Sr_4(\phi(P), \phi(Q), x) = \Sr_4(P, Q, x)(x, \alpha) = P(x, \alpha) \vspace{1mm} \\= x^4 - (\alpha + 2)\, x^3 + (2\, \alpha + 1)\, x^2 - (\alpha + 2)\, x + 2\, \alpha  $.
\end{itemize}
\end{example}

The \emph{Sturm-Habicht sequence} refines the subresultant sequence by introducing sign adjustments to form a signed subresultant sequence. Let \(P, \, Q \in \mathbb{D}[x]\) with \(p = \deg(P)\), \(q = \deg(Q)\), and \(v = p + q - 1\). The Sturm-Habicht sequence \(\{\Sth_j(P, Q)\}_{j=0,\hdots,v+1}\) is defined as in Definition 2.2 from \cite{gonzalez1998sturm}. In the case of \(Q = 1\), 
the sequence corresponds to the signed subresultant sequence of \(P\) and 
\(dP/dx\). 
The principal Sturm-Habicht coefficients \(\{\sth_j(P,Q)\}\) are the leading coefficients of \(\{\Sth_j(P,Q)\}\). Applying the \emph{generalized sign variation function} \(\mathbf{C}\) (Definition 4.1 in \cite{gonzalez1998sturm}) to \(\{\sth_j(P,1)\}\) yields the number of real roots of \(P\).

This root counting technique is central to an algorithm in \cite{bouzidi2021computation}, exploiting its specialization stability. As shown in Proposition~\ref{Prop:Res-norm}, it counts the real roots of \(n(\omega,\gamma_\star)=0\) by computing \(\{\sth_j(n,1)\}\) for \(n\) seen as a univariate in \(\omega\).

\begin{example}\label{example:strm}
We consider the bivariate polynomial $P$
studied in Example~\ref{specialization example}, where we computed the subresultant sequence of \(P\) and 
\(\partial P/\partial x\)
for \(y = \alpha\), where \(\alpha \in \mathbb{Q}\). Using the definition of the Sturm-Habicht sequence as a signed subresultant sequence, and based on the computations performed in Example~\ref{specialization example}, we obtain the principal (leading) coefficients of the sequence \(\{\sth_i(P(x, \alpha), 1)\}_{i=4,\ldots,0}=\)
$\{1, \, 4,\, 3 \, \alpha^2 - 4 \, \alpha + 4,\, \\
2\, \alpha^4 - 40\, \alpha^3 + 74\, \alpha^2 - 40\, \alpha - 128, 
-100\, (\alpha - 2)^2\, (\alpha^2 + 1)^2\}.$
By leveraging the specialization property, we can directly evaluate the principal coefficients for any value of \(y = \alpha\) without recomputing the sequence. For example, for \(\alpha = 2\), the principal coefficients become $\{\sth_i(P(x, 2), 1)\}_{ i = 4, \ldots, 0}= \{1, \, 4, \, 8, \, -200, \, 0\}$.
The sign variation \(\mathbf{C}\) applied to this sequence yields $
\mathbf{C}(1, 4, 8, -200, 0) = 1,$
indicating that \(P(x, 2)\) has exactly one real root. Similarly, for \(\alpha = 3\), the principal coefficients are
$\{\sth_i(P(x, 3), 1)\}_{, i=4, \ldots, 0}= \{1, \, 4, \, 19, \, -500, \, -1000\}.$
The sign variation \(\mathbf{C}\) applied to this sequence gives $ \mathbf{C}(1, 4, 19, -500, -1000) = 2$, indicating that \(P(x, 3)\) has two real roots.
This example highlights the utility of the specialization property: once the Sturm-Habicht sequence is computed symbolically, the principal coefficients can be evaluated for any desired value of \(\alpha\), enabling efficient real root counting for \(P(x, y=\alpha)\) without recomputing the sequence.
\end{example}

\textbf{Rational Univariate Representation (RUR).}

A key symbolic approach for solving polynomial systems with finitely many complex solutions is \emph{rational parametrization} (\cite{rouillier1999solving,rouillier2004efficient,bouzidi2015separating,bouzidi2016solving}). The \emph{Rational Univariate Representation} (RUR), introduced in \cite{rouillier1999solving} and refined for bivariate systems in \cite{bouzidi2016solving}, has recently been improved to handle demanding cases efficiently (\cite{demin2024reading}).

RUR reduces a bivariate system \(\Sigma=\{P=0,\,Q=0\}\) (with \(P, Q \in \mathbb{Z}[x,y]\)) to a univariate form by expressing solutions as roots of a polynomial \(f(t)\) and writing
\[
x = g_x(t), \quad y = g_y(t),
\]
with \(g_x(t), g_y(t) \in \mathbb{Q}(t)\).

Key steps in RUR include computing a \emph{separating linear form} (ensuring a one-to-one correspondence between univariate roots and system solutions), performing a \emph{triangular decomposition} (to hierarchically structure the system), and \emph{isolating real roots} using rational isolating boxes. This bijection guarantees that every root of \(f(t)\) (with multiplicity) corresponds to a solution of \(\Sigma\), so that solving \(f(t)=0\) recovers all solutions of $\Sigma$.

\begin{example}
Consider the polynomial \(P\) from Examples~\ref{specialization example} and \ref{example:strm}. The system
\(
\Sigma=\{P(x,y)=0,\; \partial P/\partial x=0\}
\)
defines the \(y\)-critical points of \(P(x,y)=0\). Its RUR is 
\[
x=t,\quad y=\frac{2\, (t^2-t+3)}{(3 \, t-1)(t-1)},\quad f(t)=(t-2)^2 \, (t^2+1)^2.
\]
Since \(t=2\) is the only real root of \(f\) and \((x,y)=(2,2)\) is the unique \(y\)-critical point of \(P(x,y)=0\).

\end{example}

In \cite{bouzidi2021computation}, we proposed three methods for computing the \(L_{\infty}\)-norm of a transfer matrix \(G\). The first, based on the Rational Univariate Representation (RUR), computes the maximal root \(\gamma_1\) of \(\operatorname{Lc}_\omega(n)\) (with \(n\) the polynomial associated with the \( L_{\infty} \)-norm), then solves \(\{n=0,\, \partial n/\partial\omega=0\}\) via RUR to obtain the real solutions \((\omega, \, \gamma)\). The maximum \( \gamma \) among these, \(\gamma_2\), is then compared with \(\gamma_1\), and \(\|G\|_\infty\) is set to \(\max\{\gamma_1,\gamma_2\}\). This algorithm is implemented in {\tt Maple} using the RUR code from \cite{demin2024reading} and is employed in our benchmarks.

\begin{example}
Consider the transfer function \(G\) from Example~\ref{resultants_and_norm} with $
n = 4 \, \gamma^2 \, \omega^4 + \gamma^2 \, \omega^2 + 4 \, \gamma^2 - 1$. 
The leading coefficient in \(\omega\) is \(4 \, \gamma^2\) with unique root \(\gamma = 0\). Solving \( \{n = 0, \, 
\partial n/\partial \omega 
= 0\} \)
via RUR yields \(\omega = 0\) and \(\gamma = \pm 0.5\). Hence, the maximal \(\gamma\) is \(0.5\) and \(\|G\|_\infty = 0.5\).
\end{example}

\textbf{Cylindrical Algebraic Decomposition.}
 
\emph{Cylindrical Algebraic Decomposition} (CAD) is a core method in real algebraic geometry for solving polynomial inequalities and analyzing parameter-dependent systems. It partitions \(\mathbb{R}^n\) into cells where each polynomial's sign is invariant, allowing solutions to be determined from a single representative per cell. CAD thus enables precise analysis of solution behavior in fields such as robotics, control, and optimization (see \cite{Basu} for details).

CAD construction involves two phases: projection and lifting. In projection, discriminants, resultants, and leading coefficients are computed to reduce the problem to lower dimensions while capturing critical behaviors. Lifting then reconstructs the higher-dimensional cells using sample points from these projections.
A key refinement in CAD is the use of the \emph{discriminant variety} to isolate parameters that yield non-generic behavior (e.g., infinite solutions or multiplicities). The remaining space is partitioned into cells where the number of real solutions is invariant. For instance, for \(P=x^2+a \, x+b\), the discriminant \(a^2-4 \, b=0\) (a downward parabola in \((a, \, b)\)) divides the space into regions with two roots below the curve, none above, and a double root on the curve.
The algorithm in \cite{quadrat:hal-04646145} leverages CAD to represent the maximal \(y\)-projection of \(\Sigma\)'s real solutions as a function of the parameters, enabling exact and stable parametric computation of the \(L_{\infty}\)-norm in specific regions.

\subsection{Symbolic-Numeric Approach}

Symbolic algorithms produce \emph{certified} solutions, guaranteeing correctness for any input in a finite number of steps. Unlike numerical methods, which may yield approximations or errors due to floating-point and convergence issues, symbolic methods return unambiguous results and explicitly signal failures when no solution exists.
For \(L_{\infty}\)-norm computation, symbolic methods excel in parameterized systems by partitioning the parameter space into connected regions, enabling accurate analysis even for non-convex or unbounded regions (\cite{chen2013computing,dominguez2010recent,quadrat:hal-04646145}).
Techniques such as \emph{real comprehensive triangular decomposition} and Cylindrical Algebraic Decomposition (CAD) have proven effective. For example, Chen et al. \cite{chen2013computing} extended validated numerical methods to parametric cases using border polynomials and CAD to obtain exact solutions. Our previous work (\cite{bouzidi2021computation}) addressed non-parametric \(L_{\infty}\)-norm computation via symbolic reformulation (Propositions~\ref{prop:caracterization_infty_norm} and \ref{Prop:Res-norm}) and the methods of Section~\ref{symbolic-tools}. In \cite{quadrat:hal-04646145}, we extended this approach to parameterized systems by using CAD to represent the maximal \(\gamma\)-projection of real solutions as a function of the parameters, thereby computing the \(L_{\infty}\)-norm parametrically in specified regions. These methods were implemented in \texttt{Maple} and are later used for benchmarking.

\section{Symbolic \& Numerical Methods}\label{sec:versus}

We benchmark numerical and symbolic \(L_{\infty}\)-norm computations on two cases: fixed-coefficient systems and parameterized systems, where numerical methods often struggle. In these examples, the symbolic approach is implemented in {\tt Maple} using the RUR method (\cite{bouzidi2021computation}) with the new implementation from \cite{demin2024reading}.

\begin{example}
Consider $
G(s) = \frac{s^2}{(s^2-10^{-7})(s^2-10^7)}.$
The symbolic method in {\tt Maple} computes \(\|G\|_\infty \approx 9.999998000\times10^{-8}\). In contrast, both {\tt Maple}'s {\tt NormHinf} and {\tt Matlab}'s {\tt hinfnorm} fail due to eigenvalues near zero (e.g., {\tt Matlab} drops the \(10^{-7}\) term, yielding \(\|G\|_\infty = \infty\)).
\end{example}

\begin{example}
Consider the transfer matrix
\[
G(s)=\begin{bmatrix}
\frac{10 \, (s+1)}{s^2+0.2 \, s+100} & \frac{1}{s+1} \vspace{1mm}\\
\frac{s+2}{s^2+0.1 \, s+10} & \frac{5\, (s+1)}{(s+2)(s+3)}
\end{bmatrix}.
\]
The symbolic method gives \(\|G\|_\infty=50.25\), while {\tt Maple}'s numerical method returns \(1.32\), underestimating the norm. Matlab's {\tt hinfnorm} computes \(\|G\|_\infty=50.25\) accurately. This example (Zhou and Doyle, Example 4.3) shows how sparse frequency grids in numerical methods can lead to significant errors (e.g., \(32.86\) instead of \(50.25\)).
\end{example}
\begin{example}
Consider $
G(s)=\frac{1}{(s^2+2 \, \xi s+1) \, (s+1)}.$
Adapted from Zhou and Doyle's Example 4.8, this example examines the variation of the \(L_{\infty}\)-norm with the damping ratio \(\xi\). The symbolic method computes the norm accurately for all \(\xi\); for each tenfold decrease in \(\xi\), the norm increases roughly tenfold (e.g., \(\xi=10^{-1}: 3.57\); \(10^{-2}: 35.36\); \(10^{-3}: 353.55\); \(10^{-4}: 3535.53\); \(10^{-5}: 35355.34\)). It remains robust even for \(\xi=10^{-500}\).
In contrast, {\tt Maple}'s numerical method approximates the norm for larger \(\xi\) but fails for \(\xi\le10^{-8}\), while {\tt Matlab}'s method works down to \(\xi=10^{-16}\) but fails for \(\xi<10^{-17}\).
To analyze the \(L_{\infty}\)-norm as a function of \(\xi\), we use our parameter-dependent systems algorithm (\cite{quadrat:hal-04646145}). It partitions the \(\xi\)-space into regions where the polynomial $n\in\mathbb{Q}[\xi][\omega,\gamma]$
(Definition~\ref{def:curve_n}) has a fixed number of real solutions and the roots of $
R=\operatorname{Res}\Bigl(n,\frac{\partial n}{\partial \omega},\omega\Bigr)\in\mathbb{Q}[\xi,\gamma]$ are consistently ordered. The norm is then identified by its position among these ordered roots in each region. The \(\xi\)-space is partitioned as  $$\begin{array}{lll}
C_1 = 0 < \xi < 0.5, && 
    C_2 = 0.5 < \xi < 1, \vspace{1mm}\\
    C_3 = 1 < \xi < \sqrt{5}/2, && C_4 = \xi>\sqrt{5}/2. 
\end{array}$$ 
The computed indices \(\{[C_1,7],[C_2,3],[C_3,4],[C_4,5]\}\) indicate the norm's ordered position (from bottom up) in each region. Figure~\ref{fig:root_regions_xi} 
visualizes these regions and root functions.

\begin{figure}[htbp]
\centering
\includegraphics[width=0.55\textwidth,height=0.35\textwidth,keepaspectratio]{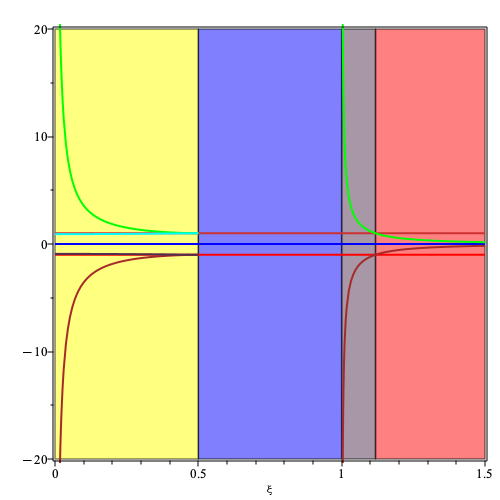}
\caption{Regions and root functions for $\xi>0$.}
\label{fig:root_regions_xi}
\end{figure}

\end{example}

\subsection{Benchmarking Methods for \( H_\infty \)-Norm Computation}\label{sec:benchmark}
This section benchmarks \(H_\infty\)-norm computation for transfer matrices using {\tt Matlab} and {\tt Maple}. We compare four methods: {\tt Matlab}’s \texttt{getPeakGain} and \texttt{hinfnorm}, {\tt Maple}’s numerical \texttt{NormHinf} from the \texttt{DynamicSystems} package, and our symbolic method based on Rational Univariate Representation (\cite{demin2024reading}). Experiments were conducted on a MacBook Pro (Apple M4 Max, 36 GB, macOS 15.1). Table~\ref{tab:benchmark_hinfinity} summarizes results for random square transfer matrices with integer coefficients, where the degree indicates the polynomial degree of both numerator and denominator. The \textbf{MatNumNorm} column shows {\tt Matlab} norms (both functions identical, all under 0.05 s), \textbf{MapNumNorm} reports {\tt Maple}’s numerical norms (all under 10 s, “FAIL” = numerical error), \textbf{MapSymNorm} gives symbolic norms with execution times ($\cdot$), and the \textbf{JuliaTime} column reports execution time of the symbolic method via PACE.jl on {\tt Julia}. The matrices and benchmark codes are available upon request. Note that the \textbf{JuliaTime} column reports only execution times; the computed \(H_\infty\)-norm values coincide exactly with those of \textbf{MapSymNorm} column.

\begin{table}[htb]
\centering
\large
\resizebox{\columnwidth}{!}{%
\begin{tabular}{lcccc}
\toprule
\textbf{Size – Degree} & \textbf{MatNumNorm} & \textbf{MapNumNorm} & \textbf{MapSymNorm} (Time (s)) & \textbf{JuliaTime} (s) \\
\midrule
2×2 – 2 & 3.0116                     & 3.0116                     & 3.0116 (0.164)   & 0.8574    \\
2×2 – 3 & 8.5769                     & 8.5769                     & 8.5769 (0.605)   & 0.1147    \\
2×2 – 4 & 8.9850                     & \textcolor{red}{8.0848}    & 8.9854 (1.399)   & 0.2525    \\
2×2 – 5 & 239.9881                   & \textcolor{red}{67.7176}   & 239.9881 (2.890) & 0.4493    \\
3×3 – 2 & 14.1745                    & FAIL                       & 14.1745 (6.503)  & 1.7485    \\
3×3 – 3 & 12.3029                    & FAIL                       & 12.3029 (32.452) & 2.0272    \\
3×3 – 4 & 97.6726                    & FAIL                       & 97.6726 (98.472) & 6.7816    \\
3×3 – 5 & 689.9393                   & FAIL                       & 689.9393 (262.320)& 11.0948   \\
4×4 – 2 & 12.5704                    & FAIL                       & 12.5705 (135.942)& 8.4516    \\
4×4 – 3 & \textcolor{red}{3.5401E+15}& FAIL                       & \textbf{50.7291} (838.939)& 42.3831  \\
4×4 – 4 & 215.6080                   & FAIL                       & 215.6080 (2739.723)& 100.8656 \\
4×4 – 5 & 253.8111                   & FAIL                       & 253.8111 (9590.923)& 355.4047 \\
5×5 – 2 & 9.5883                     & \textcolor{red}{8.5472}    & 9.5884 (1990.221)& 115.2145  \\
5×5 – 3 & 23.2985                    & \textcolor{red}{14.8156}   & 23.2985 (13485.492)& 684.6406 \\
5×5 – 4 & \textcolor{red}{6.8041E+14}& FAIL                       & \textbf{87.2824} (56303.077)& 1829.9035 \\
5×5 – 5 & \textcolor{red}{4.5986E+15}& FAIL                       & \textbf{176.1463} (164670.252)& 4239.2970 \\
\bottomrule
\end{tabular}%
}
\caption{Benchmark results for \(H_\infty\) norm computation.}
\label{tab:benchmark_hinfinity}
\end{table}

While numerical methods (\textbf{MatNumNorm} and \textbf{MapNumNorm}) scale to high degrees (up to \(n=500\)), they can yield inaccurate results (discrepancies in red) and sometimes fail. Symbolic methods are robust but computationally expensive. To improve performance, we benchmarked our symbolic approach in Julia using PACE.jl\footnote{{https://pace.gitlabpages.inria.fr/pace.jl/}}, where {\tt RationalUnivariateRepresentation.jl} handles the RUR step and optimized C routines perform isolation.


\section{Conclusion}
\label{sec:conclusion}

This paper explored symbolic computation for \(L_{\infty}\)-norm calculation in control theory, highlighting its theoretical foundation and practical impact. Comparisons with numerical techniques revealed that, although more computationally intensive, symbolic methods deliver certified accuracy and effectively complement traditional approaches. Benchmark results underscore the trade-off between efficiency and accuracy, while integration with hybrid frameworks points to promising research directions. Overall, symbolic methods are essential for advancing robust, precise control system analysis.

{\small
\bibliography{ifacconf}
}
\end{document}